\documentclass[amscd,amssymb,11pt,reqno]{amsart}
\newtheorem{Thm}{Theorem}[section]

\newtheorem{Cor}[Thm]{Corollary}
\newtheorem{Lem}[Thm]{Lemma}

\newtheorem{Def}[Thm]{Definition}
\newtheorem{Rmk}[Thm]{Remark}

\usepackage{graphicx}
\usepackage{latexsym}
\begin{document}

%\begin{Large}

\vspace{1.5 cm}

%\begin{center}{ Last modified November 7, 2004  } \end{center}

\title[Centroids and comparison of volumes.]%the Busemann-Petty problem.]
      {Centroids and comparison of volumes.}%the Busemann-Petty problem.}

\author{V.Yaskin and M.Yaskina}

\address{V.Yaskin, Department of Mathematics, University of Missouri, Columbia, MO 65211, USA}
\email{yaskinv@math.missouri.edu}

\address{M.Yaskina, Department of Mathematics, University of Missouri, Columbia, MO 65211, USA}
\email{yaskinam@math.missouri.edu}

\begin{abstract}
For $-1<p<1$ we introduce the concept of a polar $p$-centroid body ${\Gamma^*_p K}$ of a star body
$K$. We consider the question of whether ${\Gamma^*_p K}\subset {\Gamma^*_p L}$ implies
$\mathrm{vol}(L)\le \mathrm{vol}(K).$ Our results extend the studies by Lutwak in the case $p=1$
and  Grinberg, Zhang in the case $p> 1$.
\end{abstract}

\maketitle

\section{Introduction}
Let  $K$ be a star body in $\mathbb{R}^n$, then the centroid body of $K$ is a convex body $\Gamma
K$ defined by its support function:

$$h_{\Gamma K}(\xi) = \frac{1}{\mathrm{vol}(K)}\int_K |( x,\xi)| dx, \quad \xi \in \mathbb{R}^n .$$

Let  $K$ and $L$ be two origin-symmetric star bodies in $\mathbb{R}^n$ such that ${\Gamma K}\subset
{\Gamma L}$, what can be said about the volumes of $K$ and $L$? Lutwak \cite{L} proved that, if $L$
is a polar projection body then $\mathrm{vol}(K)\le \mathrm{vol}(L)$. On the other hand, if $K$ is
not a polar projection body, then there is a body $L$, so that ${\Gamma K}\subset {\Gamma L}$, but
$\mathrm{vol}(K)> \mathrm{vol}(L)$. Since in $\mathbb{R}^2$ every convex body is a polar projection
body \cite{S}, the results of Lutwak imply the following:

{\it Suppose that $K$ and $L$ are two origin-symmetric convex bodies in $\mathbb{R}^n$ such that
${\Gamma K}\subset {\Gamma L}$. If $n=2$, then we necessarily have $\mathrm{vol}(K)\le
\mathrm{vol}(L)$, while this is no longer true if $n\ge 3$.}

%In \cite{LZ} Lutwak and Zhang introduced the concept of a $p$-centroid body.
Let  $K$ be a star
body in $\mathbb{R}^n$ and $p\ge 1$, then the $p$-centroid body of $K$ is the body $\Gamma_p K$
defined by:

\begin{eqnarray}\label{eqn:p-centroid}
h_{\Gamma_p K}(\xi) = \left(\frac{1}{ \mathrm{vol}(K)}\int_K |( x,\xi) |^p dx\right)^{1/p}, \quad
\xi \in \mathbb{R}^n.
\end{eqnarray}
Note, that if $p\ge 1$, then $h_{\Gamma_p K}$ is a convex function, and, therefore, $\Gamma_p K$ is
well-defined. %(The  definition above is slightly different from that in \cite{LZ}, since they also
%have a normalization constant, which for our purposes does not matter.)
The polar of $\Gamma_p K$
is called the polar $p$-centroid body of $K$ and denoted by $\Gamma^*_p K$. Since the support
function of a body is the norm of its polar, $h_{}=\|\cdot\|_{*}$,  the polar $p$-centroid body of
$K$ is given by
\begin{eqnarray}\label{eqn:polar-p-centroid}
\|\xi\|_{\Gamma^*_p K} = \left(\frac{1}{ \mathrm{vol}(K)}\int_K |( x,\xi) |^p dx\right)^{1/p},
\quad \xi \in \mathbb{R}^n.
\end{eqnarray}

The p-centroid bodies and their polars have recently been studied by different authors, see e.g.
\cite{CG}, \cite{GZ}, \cite{L}, \cite{LYZ}, \cite{LZ}. In \cite{GZ} Grinberg and Zhang generalized
the results of Lutwak discussed in the beginning of this section. Namely, let $K$ and $L$ be two
origin-symmetric star bodies in $\mathbb{R}^n$ such that for $p\ge 1$ $${\Gamma_p K}\subset
{\Gamma_p L}.$$ They prove that if the space $(\mathbb{R}^n, \|\cdot\|_L)$ embeds in $L_p$, then we
necessarily have
$$\mathrm{vol}(K)\le \mathrm{vol}(L).$$
On the other hand, if $(\mathbb{R}^n, \|\cdot\|_K)$ does not embed in $L_p$, then there is a body
$L$ so that $\Gamma_p K \subset \Gamma_p L$, but $\mathrm{vol}(K)\le \mathrm{vol}(L).$

Note, that if $p=1$ the positive answer holds for all convex bodies in $\mathbb{R}^2$, while if
$p>1$ there is no dimension where this would always be true. The preceding remark suggests
considering $p<1$ in order to make the answer affirmative in higher dimensions.

%However this is no longer true in dimensions $n\ge 3$. One can ask a similar question about $p$th
%centroid bodies with $p> 1$. Does $\Gamma_p K \subset \Gamma_p L$ %$h_{\Gamma_p K}\le h_{\Gamma_pL}$
%imply the inequality for volumes? We prove  that the answer is negative in all dimensions.

If $p<1$, then the function $h_{\Gamma_p K}(\xi)$ in (\ref{eqn:p-centroid}) is not necessarily
convex, therefore it is not a support function, but the definition of the polar $p$-centroid body
still makes  sense, even though these bodies may be non-convex. So for all $p>-1$, $p\ne 0$ we
define the polar $p$-centroid body of a star body $K$ by the formula:
\begin{eqnarray}\label{eqn:polar-p-centroid}
\|\xi\|_{\Gamma^*_p K} = \left(\frac{1}{ \mathrm{vol}(K)}\int_K |( x,\xi) |^p dx\right)^{1/p},
\quad \xi \in \mathbb{R}^n.
\end{eqnarray}
For $p=0$, this definition looks as follows (if we send $p\to 0$):
\begin{eqnarray}\label{eqn:polar-0-centroid}
\|\xi\|_{\Gamma^*_0 K} = \exp\left(\frac{1}{ \mathrm{vol}(K)}\int_K \ln |( x,\xi) | dx\right),
\quad \xi \in \mathbb{R}^n.
\end{eqnarray}

Now we can ask the question discussed above for  all $p>-1$.  Namely, suppose that
\begin{eqnarray}\label{condition}
{\Gamma^*_p L}\subset {\Gamma^*_p K},
\end{eqnarray} %What can be said about the volumes of $K$ and $L$?
for origin-symmetric star bodies $K$ and $L$. Does it follow that we have an inequality for the
volumes of $K$ and $L$? In this paper we show that if $(\mathbb{R}^n, \|\cdot\|_L)$ embeds in
$L_p$, $p>-1$, then we have $\mathrm{vol}(K)\le \mathrm{vol}(L)$. However if $(\mathbb{R}^n,
\|\cdot\|_K)$ does not embed in $L_p$, we construct counterexamples to the latter result.

These results can also be reformulated as follows: \newline (i) If $0<p<1$, then in $\mathbb{R}^2$
the condition (\ref{condition}) implies that $\mathrm{vol}(K)\le \mathrm{vol}(L)$, while this is no
longer true in dimensions $n\ge 3$. \newline (ii) If $-1<p\le 0$,  (\ref{condition}) implies that
$\mathrm{vol}(K)\le \mathrm{vol}(L)$ if and only if  $n\le 3$.

%Finally, ${\Gamma^*_0 L}\subset {\Gamma^*_0 K}$ implies $\mathrm{vol}(K)\ge \mathrm{vol}(L)$ in
%$\mathbb{R}^3$ and not higher.

Clearly the integral in (\ref{eqn:polar-p-centroid}) diverges if $p\le -1$, but still we can make
sense of this integral considering fractional derivatives. Indeed, if $-1<p<0$
\begin{eqnarray*}\frac{1}{\mathrm{vol}(K)} \int_K |( x,\xi) |^{p} dx
&=& \frac{1}{\mathrm{vol}(K)}\int_{-\infty}^{\infty} |z|^p \int_{(x,\xi)=z} \chi({\|x\|_K}) dx \ dz\\
&=&\frac{1}{\mathrm{vol}(K)}\int_{-\infty}^{\infty} |z|^p A_{K,\xi}(z) dz \\&=&
\frac{2\Gamma(p+1)}{\mathrm{vol}(K)}  A^{(-p-1)}_{K,\xi}(0),
\end{eqnarray*}
where $A_{K,\xi}(z)$ is  the parallel section function of $K$, and $A^{(-p-1)}_{K,\xi}(0)$ is its
fractional derivative at zero. (For details on fractional derivatives, see e.g. \cite[Section
2.6]{K3}). So, in such terms our problem can be written as follows:

Suppose $K$ and $L$ are  two origin-symmetric star bodies, so that for all $\xi \in S^{n-1}$:
$$\frac{A^{(-p-1)}_{K,\xi}(0)}{\mathrm{vol}(K)}\le \frac{A^{(-p-1)}_{L,\xi}(0)}{\mathrm{vol}(L)}.$$
Do we necessarily have an inequality for the volumes of $K$ and $L$?

Note that Koldobsky already considered such inequalities (see e.g. \cite{K1}) without dividing by
volumes. So, for $-1<p<0$ the positive part of our results can also be obtained from the results of
Koldobsky, but we give our own proof. The case $p=-1$   leads to the following modification of the
Busemann-Petty problem. Let $K$ and $L$ be two convex origin-symmetric bodies in $\mathbb{R}^n$
such that

$$\frac{\mathrm{vol}_{n-1} (K\cap
\xi^\perp)}{\mathrm{vol}(K)}\le \frac{\mathrm{vol}_{n-1} (L\cap \xi^\perp)}{\mathrm{vol}(L)}.$$
Does this imply an inequality for the volumes of $K$ and $L$?

It is easy to show that in dimensions $n\le 4$ we have $\mathrm{vol}(L) \le \mathrm{vol}(K).$ The
proof is almost identical  to that of the original Busemann-Petty problem from \cite{GKS}. The
counterexamples in dimensions $n\ge 5$ from \cite{GKS} also  work in this situation.

In view of all these remarks one can consider our results as a certain bridge between the results
of Lutwak-Grinberg-Zhang about $p$-centroid bodies and the results of Busemann-Petty type obtained
by Koldobsky.

%\section{Preliminaries}

\section{ Centroid inequalities for $-1<p<1$, $p\ne 0$.}

The {\it Minkowski functional} of a star-shaped origin-symmetric body $K\subset \mathbb R^n$ is
defined as
$$\|x\|_K=\min \{a\ge 0: x \in aK \}.$$
We denote by $(\mathbb{R}^n,\|\cdot\|_K)$ the Euclidean space equipped with the Minkowski
functional of the body $K.$ Clearly, $(\mathbb{R}^n,\|\cdot\|_K)$ is a normed space if and only if
the body $K$ is convex.

The {\it support function} of a convex body $K$ in $\mathbb{R}^n$ is defined by
$$h_K(x)=\max_{\xi \in K}(x,\xi),\ \ x\in \mathbb{R}^n. $$
If $K$ is origin-symmetric, then $h_K$ is the Minkowski norm of the polar body $K^*$.

A well-known result going back to P.L\'{e}vy, (see \cite[p. 189]{BL} or \cite[Section 6.1]{K3}), is
that a space $(\mathbb{R}^n, \|\cdot\|)$ embeds into $L_p$, $p>0$ if and only if there exists a
finite Borel measure $\mu$ on the unit sphere so that, for every $x\in\mathbb{R}^n$,
\begin{eqnarray}{\label{Def:L_p>0}}
\|x\|^p=\int_{S^{n-1}} |(x, \xi ) |^p d\mu(\xi).
\end{eqnarray}
On the other hand, this can be considered as the definition of embedding in $L_p$, $-1<p<0$ (cf.
\cite{K5}).

%We say that $(\mathbb{R}^n, \|\cdot\|)$ embeds into $L_p$, $-1<p<0$ if there exists a finite Borel
%measure $\mu$ on the unit sphere so that, for every $x\in\mathbb{R}^n$,
%$$\|x\|^p=\int_{S^{n-1}} |(x, \xi ) |^p d\mu(\xi).$$

It was proved in \cite{K4} that a space $(\mathbb{R}^n,\|\cdot\|)$ embeds isometrically in $L_p,\
p>0,$ $p \notin 2\mathbb{N}$ if and only if the Fourier transform of the function
$\Gamma(-p/2)\|x\|^{p}$ (in the sense of distributions) is a positive distribution outside of the
origin. If $-n<p<0$ a similar fact was proved in \cite{K5}: a space $(\mathbb{R}^n,\|\cdot\|)$
embeds in $L_p$ if and only if the Fourier transform of $\|\cdot\|^p$ is a positive distribution in
the whole $\mathbb{R}^n.$

% Throughout the paper, we write $(\mathbb{R}^n,\|\cdot\|)$ meaning that
%$\|\cdot\|$ is the Minkowski functional of some origin-symmetric star body.

%Also we will need the following
%\begin{Lem}\cite[Section 2.2]{MP}
%Let $K$ and $M$ be two origin-symmetric convex bodies in $\mathbb{R}^n$ and let $p>0$. Then the
%following inequalities hold:

%\begin{eqnarray}\label{eqn:MP} \left(\frac{1}{\mathrm{vol}(K)}\int_K ||x||^p_M dx\right)^{1/p}\ge
%\left(\frac{n}{n+p}\right)^{1/p} \left(\frac{\mathrm{vol}(K)}{\mathrm{vol}(M)}\right)^{1/n},
%\end{eqnarray}
%and
%\begin{eqnarray}\label{eqn:logMP} \exp\left(\frac{1}{\mathrm{vol}(K)}\int_K \ln ||x||_M dx\right)\ge
%e^{-1/n} \left(\frac{\mathrm{vol}(K)}{\mathrm{vol}(M)}\right)^{1/n}.
%\end{eqnarray}
%\end{Lem}

 %In this section we prove such results for $-1<p< 1$, $p\neq 0$.

Now we are ready to prove our first result.

\begin{Thm}\label{Thm:2D}
Let $-1<p< 1$, $p\neq 0$. Let $K$ and $L$ be  origin-symmetric convex bodies in $\mathbb{R}^n$, so
that $(\mathbb{R}^n, \|\cdot\|_K)$ embeds in $L_p$ and
\begin{eqnarray}\label{condition1}
{\Gamma^*_p K}\subset {\Gamma^*_p L}.
\end{eqnarray}
 Then $\mathrm{vol}(L)\le \mathrm{vol}(K).$
\end{Thm}

\proof First let us prove the case  $0<p<1$. Since $(\mathbb{R}^n, \|\cdot\|_K)$ embeds in $L_p$,
there exists a measure $\mu_K$ on the unit sphere $S^{n-1}$ such that
$$\|x\|_K^p=\int_{S^{n-1}}|( x,\xi)|^pd \mu_K(\xi). $$
%this is a well-known result of P.L\'{e}vy, see \cite[p. 189]{BL} or \cite[Section 6.1]{K3}

Note that (\ref{condition1})  can be written as
\begin{eqnarray}
\frac{1}{\mathrm{vol}(L)}\int_L |( x,\xi)|^pd x\le \frac{1}{\mathrm{vol}(K)}\int_K |( x,\xi)|^pd x,
\end{eqnarray}

Integrating both sides of the last inequality over $S^{n-1}$ with the measure $\mu_K$, we get
$$\frac{1}{\mathrm{vol}(L)} \int_{S^{n-1}}\int_L |( x,\xi)|^p d x\ d \mu_K(\xi)\le
 \frac{1}{\mathrm{vol}(K)}\int_{S^{n-1}}\int_K |( x,\xi)|^p d x\ d \mu_K(\xi).$$

Applying Fubini's Theorem,
\begin{equation}\label{Thm:1part}
\frac{1}{\mathrm{vol}(L)} \int_L \|x\|_K^p d x\le  \frac{1}{\mathrm{vol}(K)} \int_K \|x\|_K^p d x.
\end{equation}

Note that
\begin{eqnarray*}
\int_K \|x\|_K^p d x &=&\int_{S^{n-1}}\left(\int_0^{\|\theta\|_K^{-1}}  \| r\theta\|_K^p \ r^{n-1}
d r
\right) d\theta\\
&=&\frac{1}{n+p}\int_{S^{n-1}}\|\theta\|^{-n}_K d \theta=\frac{n}{n+p}\mathrm{vol}(K).
\end{eqnarray*}

Therefore, (\ref{Thm:1part}) can be rewritten as
\begin{eqnarray*}
\frac{1}{\mathrm{vol}(L)} \int_L \|x\|_K^p d x \le \frac{n}{n+p}.
\end{eqnarray*}

%\begin{eqnarray*}
%\int_L ||x||_Kd x &=&\int_{S^1}\left(\int_0^\infty r^2 ||\theta||_K \chi(r ||\theta||_L) d r
%\right) d\theta\\
%&=&\frac{1}{3}\int_{S^1}||\theta||^{-3}_L||\theta||_K d \theta\\
%\end{eqnarray*}

Using the inequality
\begin{eqnarray}\label{eqn:MP}
\frac{1}{\mathrm{vol}(L)}\int_L \|x\|^p_K dx\ge  \frac{n}{n+p}
\left(\frac{\mathrm{vol}(L)}{\mathrm{vol}(K)}\right)^{p/n}
\end{eqnarray}
from \cite[Section 2.2]{MP}, we get
\begin{eqnarray*}
\frac{n}{n+p} &\ge& \frac{1}{\mathrm{vol}(L)} \int_L \|x\|_K^p d x \ge \frac{n}{n+p}
\left(\frac{\mathrm{vol}(L)}{\mathrm{vol}(K)}\right)^{p/n},
\end{eqnarray*}
therefore $\mathrm{vol}(L)\le \mathrm{vol}(K),$ which proves the theorem for $0<p<1$.

Now consider $-1<p<0$. In this case (\ref{condition1}) is equivalent to
\begin{eqnarray}\label{cond:p<0}
\frac{1}{\mathrm{vol}(L)}\int_L |( x,\xi)|^pd x\ge \frac{1}{\mathrm{vol}(K)}\int_K |( x,\xi)|^pd x,
\end{eqnarray}

Since $(\mathbb{R}^n, \|\cdot\|_K)$ embeds into $L_p$, $p>-1$, there exists a measure $\mu_K$ on
the unit sphere such that
$$\|x\|_K^p=\int_{S^{n-1}}|( x,\xi)|^pd \mu_K(\xi). $$

Integrating both sides of (\ref{cond:p<0}) over $S^{n-1}$ with the measure $\mu_K$ and using the
same argument as in the first part of the proof, we get
\begin{equation}\label{eneq:M-K}
 \frac{1}{\mathrm{vol}(L)}\int_L \| x\|^p_K d x \ge \frac{n}{n+p} .
\end{equation}

Passing to spherical coordinates and applying H\"{o}lder's inequality
\begin{eqnarray*}
\int_L \| x\|^p_K d x &=&\int_{S^{n-1}}\left(\int_0^{\|\theta\|^{-1}_L}r^{n+p-1}\|\theta\|_K^p dr
\right) d\theta\\
&=&\frac{1}{n+p}\int_{S^{n-1}}\|\theta\|^{-n-p}_L \|\theta\|_K^p  d\theta\\
&\le&\frac{1}{n+p}\left(\int_{S^{n-1}}\|\theta\|^{-n}_L
d\theta\right)^{(n+p)/n}\left(\int_{S^{n-1}}
\|\theta\|_K^{-n} d\theta\right)^{-p/n}\\
&=&\frac{n}{n+p}\left(\mathrm{vol}(L)\right)^{(n+p)/n}\left(\mathrm{vol}(K)\right)^{-p/n}.
\end{eqnarray*}

So (\ref{eneq:M-K}) can be written as
\begin{eqnarray*}
1&\le&\frac{1}{\mathrm{vol}(L)}\left(\mathrm{vol}(L)\right)^{(n+p)/n}\left(\mathrm{vol}(K)\right)^{-p/n}\\
&=&\left(\mathrm{vol}(L)\right)^{p/n}\left(\mathrm{vol}(K)\right)^{-p/n}.
\end{eqnarray*}

Therefore, using the fact that $p<0$, we get $\mathrm{vol}(L)\le \mathrm{vol}(K)$.

\qed

Since all 2-dimensional spaces embed in $L_1$, and therefore in $L_p$ with $-2<p<1$ (see e.g.
\cite[Chapter 6]{K3}), and all 3-dimensional spaces embed in $L_0$, and therefore in $L_p$ with
$-3<p<0$ (see \cite{KKYY}), we have the following
\begin{Cor}
Let $K$ and $L$ be origin-symmetric convex bodies in $\mathbb{R}^n$, so that ${\Gamma^*_p K}\subset
{\Gamma^*_p L}.$ Then

i) if $0<p<1$, we necessarily have $\mathrm{vol}(L)\le \mathrm{vol}(K)$ in dimension $n=2$,

ii) if $-1<p<0$, we necessarily have $\mathrm{vol}(L)\le \mathrm{vol}(K)$ in dimensions $n=2$ and
$3$.
\end{Cor}

In order to show a negative counterpart of Theorem \ref{Thm:2D}, we need some lemmas. The following
Lemma is \cite[Corollary 3.15]{K3} with $k=0$ and $p=-q-1$.

\begin{Lem}\label{Lem:homog}
Let $-1<p<1$, $p \neq 0$. For an origin-symmetric convex body $K$ in $\mathbb{R}^n$ we have
%$$\left(||x||_K^{-n-1}\right)^\wedge (\xi)= -\pi(n+1)A^{(-2)}_{K,\xi}(0)$$
%$$\left(\|x\|_K^{-n-p}\right)^\wedge (\xi)= \Gamma(-p)\sin\left(\frac{\pi(p+1)}{2}\right)\int_{K}|( x,\xi)|^p dx.$$
%$$\left(||x||_K^{-n-1}\right)^\wedge (\xi)= -\pi(n+1)\int_{K}|( x,\xi)| dx. $$
\begin{eqnarray*}
\left(\|x\|_K^{-n-p}\right)^\wedge (\xi)= -\frac{\pi}{2\Gamma(p+1)\sin\left({\pi
p}/{2}\right)}\int_{S^{n-1}}|( \theta,\xi)|^p \ \|\theta\|_K^{-n-p} d\theta.
\end{eqnarray*}
\end{Lem}
%\proof
We will use this formula in the following form:
\begin{eqnarray*}
\left(\|x\|_K^{-n-p}\right)^\wedge (\xi)= -\frac{\pi(n+p)}{2\Gamma(p+1)\sin\left({\pi
p}/{2}\right)}\int_{K}|( x,\xi)|^p dx.
\end{eqnarray*}

Also we can write this formula in terms of fractional derivatives of the parallel section function
of $K$. Recall that the parallel section function of a an origin-symmetric star body $K$ is defined
by
$$A_{K,\xi}(z)=\displaystyle\int_{( x, \xi)=z}\chi(\|x\|_K)dx.$$ For $-1<q<0$ the fractional
derivative of this function at zero is defined by
\begin{eqnarray*}
A^{(q)}_{K,\xi}(0)= \frac{1}{2\Gamma(-q)}\int_{-\infty}^{\infty}|z|^{-1-q} A_{K,\xi}(z) dz
                 =\frac{1}{2\Gamma(-q)}\int_{K}|( x,\xi)|^{-1-q} dx
\end{eqnarray*}
In fact one can see that this can be analytically extended to $q<-1$. Therefore Lemma
\ref{Lem:homog} can be reformulated as follows. Let $-1<p<1$, $p \neq 0$, then
$$\left(\|x\|_K^{-n-p}\right)^\wedge (\xi)= -\frac{\pi(n+p)}{\sin(\pi p/2)}A^{(-p-1)}_{K,\xi}(0).$$

Note, that for $-1<p<0$ this formula was proved in \cite{GKS}.

Now recall a version of Parseval's formula on the sphere proved by Koldobsky \cite{K2}.

\begin{Lem}\label{Lem:Parseval}
If $K$ and $L$ are origin-symmetric infinitely smooth bodies in $\mathbb{R}^n$ and $0<p<n$, then
$(\|x\|_K^{-p})^\wedge$ and $(\|x\|_L^{-n+p})^\wedge$ are continuous functions on $S^{n-1}$ and
$$\int_{S^{n-1}} \left(\|x\|_K^{-p}\right)^\wedge (\xi) \left(\|x\|_L^{-n+p}\right)^\wedge
(\xi)d\xi= (2\pi)^n \int_{S^{n-1}} \|x\|_K^{-p}\|x\|_L^{-n+p} dx.$$
\end{Lem}
\begin{Rmk}\label{Rmk:Parseval}{\rm
A proof of this formula via spherical harmonics was given in \cite{K1}. Repeating this proof word
by word and using the above definition of the fractional derivative of order $q<-1$, one can easily
extend this result to $-1<p<0$.}
\end{Rmk}
%Using the Parseval's formula we can  construct counterexamples to Theorem \ref{Thm:2D} in
%$\mathbb{R}^3$

Now we prove a negative counterpart of Theorem \ref{Thm:2D}.
\begin{Thm}\label{p-counter}
Let $L$ be an infinitely smooth origin-symmetric strictly convex body in $\mathbb{R}^n$, for which
$(\mathbb{R}^n,\|\cdot\|_L)$ does not embed in $L_p$, $-1<p<1$, $p\ne 0$. Then there exists an
origin-symmetric convex body $K$ in $\mathbb{R}^n$ such that
\begin{eqnarray*}%\label{condition}
{\Gamma^*_p K}\subset {\Gamma^*_p L}.
\end{eqnarray*}
but $$\mathrm{vol}(L)> \mathrm{vol}(K).$$

\end{Thm}
\proof %Let $K$ be the unit ball of the space $l^3_q$, $q>2$. It is known that $l^3_q$ with $q>2$
%does not embed in $L_1$, therefore $||x||_K$ is not positive definite.

%Apply the usual argument to construct another body $L$...

First consider $0<p<1$. Since $(\mathbb{R}^n,\|\cdot\|_L)$ does not embed in $L_p$, there exists a
$\xi \in S^{n-1}$ such that $\left(\|x\|^p_L\right)^\wedge (\xi)$ is positive, for more details see
\cite{K4}. Because $\left(\|x\|^p_L\right)^\wedge (\theta)$ is a continuous function on $S^{n-1}$,
there exists a neighborhood of $\xi$ where it is positive. Define
$$\Omega=\{ \theta \in S^{n-1}: \left(\|x\|^p_L\right)^\wedge (\theta)>0\}.$$
Choose a non-positive infinitely-smooth even function $v$ supported on $\Omega$. Extend $v$ to a
homogeneous function $|x|_2^{-n-p}v(x/|x|_2)$ of degree $-n-p$ on $\mathbb{R}^n$. By  \cite[Chapter
3]{K3},  the Fourier transform of $ |x|_2^{-n-p}v(x/|x|_2)$ is equal to $|x|_2^p \ g(x/|x|_2)$ for
some infinitely smooth function $g$ on $S^{n-1}$.

Define a body $K$ by
$$\|x\|^{-n-p}_K = \|x\|^{-n-p}_L+\epsilon |x|_2^{-n-p} g(x/|x|_2)$$
for some small $\epsilon$ so that the body $K$ is convex (see e.g. the perturbation argument from
\cite[Section 5.1]{K3}). Applying the Fourier transform to both sides we get

$$\left(\|x\|^{-n-p}_K\right)^\wedge(\xi) = \left(\|x\|^{-n-p}_L\right)^\wedge(\xi)+\epsilon (2\pi)^n |\xi|_2^p v(\xi/|\xi|_2). $$

So using the formula from Lemma \ref{Lem:homog}
$$\left(\|x\|_K^{-n-p}\right)^\wedge (\xi)= \Gamma(-p)\sin\left(\frac{\pi(p+1)}{2}\right)\int_{K}|( x,\xi)|^p
dx$$

%$$\left(||x||^{-n-1}_M\right)^\wedge(\xi)=\frac{-\pi(n+1)}{2}\int_{K}|(
%x, \xi)| dx$$
we have \begin{equation}\label{eqn:ineq} \int_{L}|( x, \xi)|^p dx<\int_{K}|( x, \xi)|^p dx
.\end{equation}

Consider the integral
\begin{eqnarray}\label{3example:RHS}
&&\hspace{-1cm}\int_{S^{n-1}}\left(\|x\|^p_L\right)^\wedge (\xi)\left(\|x\|_K^{-n-p}\right)^\wedge (\xi)d\xi \nonumber\\
&=&\int_{S^{n-1}} \left( \|x\|^p_L\right)^\wedge (\xi) \left(\|x\|_L^{-n-p}\right)^\wedge (\xi)d\xi + \epsilon(2\pi)^n \int_{S^{n-1}}\left(\|x\|^p_L\right)^\wedge (\xi) v(\xi)d\xi\nonumber \\
&<&\int_{S^{n-1}} \left( \|x\|^p_L\right)^\wedge (\xi) \left(\|x\|_L^{-n-p}\right)^\wedge (\xi)d\xi \nonumber \\
&=&(2\pi)^n \int_{S^{n-1}} \|x\|^p_L \|x\|_L^{-n-p}dx= (2\pi)^n n {\mathrm{vol}(L)}.
\end{eqnarray}
Here we used a version of Parseval's formula (Lemma \ref{Lem:Parseval} and Remark
\ref{Rmk:Parseval}) and the fact that $v$ is negative on $\Omega$.

On the other hand, again using Parseval's formula and (\ref{eqn:MP})
\begin{eqnarray}\label{3example:LHS}
&&\int_{S^{n-1}}\left(\|x\|^p_L\right)^\wedge (\xi)\left(\|x\|_K^{-n-p}\right)^\wedge (\xi)d\xi=
(2\pi)^n \int_{S^{n-1}}\|x\|^p_L \|x\|_K^{-n-p} d x \nonumber \\
&&=(2\pi)^n(n+p)\int_K \|x\|^p_L d x \ge (2\pi)^n n {\mathrm{vol}}(K)\left(\frac{{\mathrm{vol}}(L)}
{{\mathrm{vol}}(L)}\right)^{p/n}
\end{eqnarray}
Combining (\ref{3example:RHS}) and (\ref{3example:LHS}) we get
\begin{eqnarray}\label{eqn:volumes}
{\mathrm{vol}}(K) <{\mathrm{vol}}(L).
\end{eqnarray}
Now from (\ref{eqn:volumes})
and (\ref{eqn:ineq}) %$\int_{K}|( x, \xi)|^p dx<\int_{M}|( x, \xi)|^p dx$, the inequality we have already obtained,
it follows that
$$\displaystyle\frac{1}{\mathrm{vol}(L)}\int_L |(
x,\xi)|^pd x\le \frac{1}{\mathrm{vol}(K)}\int_K |( x,\xi)|^pd x ,$$ which is equivalent to
\begin{eqnarray*}
{\Gamma^*_p K}\subset {\Gamma^*_p L}.
\end{eqnarray*}

Now consider the case $-1<p<0$. Since $(\mathbb{R}^n,\|\cdot\|_L)$ does not embed in $L_p$,  there
exists a $\xi \in S^{n-1}$ such that $\left(\|x\|^p_L\right)^\wedge (\xi)$ is negative, see
\cite{K4}. Define $$\Omega=\{ \theta \in S^{n-1}: \left(\|x\|^p_L\right)^\wedge (\theta)<0\}$$ and
choose $v(\theta)$ the same way as in the first part.

Define a body $K$ by
$$\frac{\|x\|^{-n-p}_K}{\mathrm{vol}(K)} =\frac{\|x\|^{-n-p}_L}{\mathrm{vol}(L)}+\epsilon |x|_2^{-n-p} g(x/|x|_2)$$
for some small $\epsilon$ so that the body $K$ is convex. Applying Fourier transform to both sides
we get

$$\frac{1}{\mathrm{vol}(K)}\left(\|x\|^{-n-p}_K\right)^\wedge(\xi) = \frac{1}{\mathrm{vol}(L)}\left(\|x\|^{-n-p}_L\right)^\wedge(\xi)
+\epsilon (2\pi)^n |\xi|_2^p v(\xi/|\xi|_2). $$

Again using the formula from  Lemma \ref{Lem:homog} and the fact that $v(\theta)$ is non-positive,
we have
$$\frac{1}{\mathrm{vol}(K)}\int_{K}|( x, \xi)|^p dx<\frac{1}{\mathrm{vol}(L)}\int_{L}|( x, \xi)|^p dx ,$$
which is the same as
\begin{eqnarray*} {\Gamma^*_p K}\subset {\Gamma^*_p L},
\end{eqnarray*} since $-1<p<0$.

Consider the integral

%\begin{eqnarray}\label{Nexample:RHS}
%&&
$$\hspace{-2cm}\frac{1}{\mathrm{vol}(K)}\int_{S^{n-1}}\left(\|x\|^p_L\right)^\wedge
(\xi)\left(\|x\|_K^{-n-p}\right)^\wedge (\xi)d\xi$$ %=\nonumber\\
$$=
\frac{1}{\mathrm{vol}(L)}\int_{S^{n-1}} \left( \|x\|^p_L\right)^\wedge (\xi) \left(\|x\|_L^{-n-p}\right)^\wedge (\xi)d\xi \nonumber \\
+ \epsilon (2\pi)^n\int_{S^{n-1}}\left(\|x\|^p_L\right)^\wedge (\xi) v(\xi)d\xi$$ %\nonumber \\
\begin{equation}\label{Nexample:RHS}
>\frac{1}{\mathrm{vol}(L)}\int_{S^{n-1}} \left( \|x\|^p_L\right)^\wedge (\xi) \left(\|x\|_L^{-n-p}\right)^\wedge (\xi)d\xi %\nonumber \\
=(2\pi)^n n .
\end{equation}
%\end{eqnarray}
Here we used  Parseval's formula and the fact that $v$ is negative on $\Omega$.

On the other hand, again using Parseval's formula and H\"{o}lder's inequality
\begin{eqnarray}\label{Nexample:LHS}
&&\hspace{-2cm}\int_{S^{n-1}}\left(\|x\|^p_L\right)^\wedge (\xi)\left(\|x\|_K^{-n-p}\right)^\wedge
(\xi)d\xi=
(2\pi)^n \int_{S^{n-1}}\|x\|^p_L \|x\|_K^{-n-p} d x \nonumber \\
&\le& (2\pi)^n \left(\int_{S^{n-1}}\|x\|^{-n}_L dx\right)^{-p/n}\left(\int_{S^{n-1}} \|x\|_K^{-n} d x\right)^{(n+p)/n} \nonumber \\
&=& (2\pi)^n n\left(\mathrm{vol}(L)\right)^{-p/n}\left({\mathrm{vol}(K)}\right)^{(n+p)/n}.
\end{eqnarray}

So combining (\ref{Nexample:RHS}) and (\ref{Nexample:LHS}) we get $\mathrm{vol}(L)>
\mathrm{vol}(K)$.

 \qed

\begin{Cor} The result of Theorem \ref{p-counter} can be formulated as follows:

i) Let $-1<p<0$. There exist origin-symmetric convex bodies $K$ and $L$ in $\mathbb{R}^4$, so that
${\Gamma^*_p K}\subset {\Gamma^*_p L},$ but $\mathrm{vol}(L)> \mathrm{vol}(K).$

ii) Let $0<p<1$. There exist origin-symmetric convex bodies $K$ and $L$ in $\mathbb{R}^3$, so that
${\Gamma^*_p K}\subset {\Gamma^*_p L},$ but $\mathrm{vol}(L)> \mathrm{vol}(K).$
\end{Cor}
\noindent{\bf Proof.} Consider only the case $-1<p<0$, the other case is similar. In view of the
previous theorem it is enough to construct an origin-symmetric infinitely smooth convex body  $L\in
\mathbb{R}^4$ for which the distribution $( \|x\|_L^{p})^\wedge$ is not positive. The construction
will be similar to that from \cite{GKS}.

Define $f_N(x)=(1-x^2-Nx^4)^{1/3}$, let $a_N>0$ be such that $f_N(a_N)=0$ and $f_N(x)>0$ on the
interval $(0,a_N)$. Define a body $L$ in $\mathbb{R}^4$ by
$$L=\{(x_1,x_2,x_3,x_4)\in \mathbb{R}^4: \, x_4\in [-a_N,a_N]\,  \mathrm{ and }\,
\sqrt{x_1^2+x_2^2+x_3^2}\le f_N(x_4) \}.$$ The body $L$ is strictly convex and infinitely smooth.

By the formula
$$A_{L,\xi}^{(q)}(0)= \frac{\cos{\frac{\pi q}{2}}}{\pi (n-q-1)} \left(
\|x\|_L^{-n+q+1} \right)^\wedge (\xi)$$ from \cite{GKS} and the definition of fractional
derivatives, we get
\begin{eqnarray*}\left(
\|x\|_L^p \right)^\wedge (\xi) &=& \frac{\pi p}{\cos{\frac{\pi(3+p)}{2}}}A_{L,\xi}^{(3+p)}(0)\\
&=&\frac{\pi
p}{\Gamma(-3-p)\cos{\frac{\pi(3+p)}{2}}}\int_0^{\infty}\frac{A_{L,\xi}(z)-A_{L,\xi}(0)-
A''_{L,\xi}(0)\frac{z^2}{2}}{z^{4+p}}dz.
\end{eqnarray*}
Note that the coefficient in the latter formula is positive, therefore it is enough to show that
the integral is negative.

 The function $A_{L,\xi}$ can easily be computed: $$A_{L,\xi}(x)=\frac{4\pi}{3} (1-x^2-Nx^4).$$
We have
\begin{eqnarray*}&&\int_0^{\infty}\frac{A_\xi(z)-A_\xi(0)-
A''_\xi(0)\frac{z^2}{2}}{z^{4+p}}dz=\\
&&=\frac{4\pi}{3}\left(-\frac{1}{1+p}
Na_N^{1+p}+\frac{1}{(1+p)a_N^{(1+p)}}-\frac{1}{(3+p)a_N^{3+p}}\right).
\end{eqnarray*}
The latter is  negative for $N$ large enough, because $N^{1/4}\cdot a_N\to 1$ as $N\to \infty$.

\qed

\section{ Centroid  inequalities for $p=0$.}

In this section we extend the results of the previous section to $p=0$. First we need some
preliminary results.  The concept of embedding in $L_0$ was introduced in \cite{KKYY}:

\begin{Def}\label{Def:Lzero}We say that a space $(\mathbb{R}^n, \|\cdot\|)$ embeds in $L_0$ if there exist a finite Borel
measure $\mu$ on the sphere $S^{n-1}$ and a constant $C\in \mathbb{R}$ so that, for every
$x\in\mathbb{R}^n$,
\begin{equation} \label{logrepr}
\ln \|x\| =\int_{S^{n-1}} \ln |(x, \xi ) | d\mu(\xi) + C.
\end{equation}
\end{Def}

It follows directly from the definition that $\mu$ is a probability measure, and the constant $C$
equals
\begin{eqnarray}\label{const}
C= \frac{1}{|S^{n-1}|}\int_{S^{n-1}} \ln\|x\| dx
-\frac{1}{2\sqrt{\pi}}\Gamma'(1/2)+ \frac{1}{2}\frac{\Gamma'(n/2)}{\Gamma(n/2)}.
\end{eqnarray}

Also it was proved that if $K$ is an infinitely smooth body then $\left(\ln
\|x\|_K\right)^\wedge(\xi)$ is a  homogeneous of degree $-n$ function on $\mathbb{R}^n\setminus
\{0\}$, as seen from the following

\begin{Thm}\label{Thm:n-1-deriv} \cite[Theorem 4.1]{KKYY}
Let $K$ be an infinitely smooth origin-symmetric star body in $\mathbb{R}^n$. Extend
$A_{K,\xi}^{(n-1)}(0)$ to a homogeneous function of degree $-n$ of the variable $\xi \in
\mathbb{R}^n \setminus \{0\}$. Then

i) if $n$ is odd

$$\left( \ln \|x\|_K\right)^\wedge (\xi)=(-1)^{(n+1)/2}\pi A_{K,\xi}^{(n-1)}(0),  \ \ \xi \in \mathbb{R}^n\setminus \{0\}$$

ii) if $n$ is even, then for    $\xi \in \mathbb{R}^n\setminus \{0\}$,

$$\left( \ln \|x\|_K\right)^\wedge (\xi)=a_n\int_0^{\infty}\frac{A_\xi(z)-A_\xi(0)-
A''_\xi(0)\frac{z^2}{2}-...-A^{n-2}_\xi(z)\frac{z^{n-2}}{(n-2)!}}{z^n}dz,$$ where
$a_n=2(-1)^{n/2+1}(n-1)!$

\end{Thm}

%Homogeneity

%\begin{eqnarray}\label{eqn:homog}
%\langle \mu(x),\phi(x/t)\rangle &=& - \langle \left(\ln\|x\|\right)^\wedge(x),\phi(x/t)\rangle \nonumber\\
%&=& -\int_{\mathbb{R}^n}\ln\|z\|\hat{\phi}(t z)t^n dz \nonumber\\
%&=& -\int_{\mathbb{R}^n}\hat{\phi}(\tilde{x})\ln\|\frac{1}{t}\tilde{x}\|d\tilde{x} \nonumber\\
%&=& -\int_{\mathbb{R}^n}\hat{\phi}(\tilde{x})\ln\|\tilde{x}\|d\tilde{x}+
%\ln|t|\int_{\mathbb{R}^n}\hat{\phi}(\tilde{x})d\tilde{x} \nonumber \\
%&=& -\int_{\mathbb{R}^n}\hat{\phi}(\tilde{x})\ln\|\tilde{x}\|d\tilde{x} \nonumber \\
%&=&\langle \mu(x), \phi(x)\rangle.
%\end{eqnarray}

In particular, for an infinitely smooth origin-symmetric star body $K$,  $\left(\ln
\|x\|_K\right)^\wedge(\xi)$ is a continuous function on $S^{n-1}$,  and moreover the measure in
Definition \ref{Def:Lzero} equals

$$d\mu(\xi)=-\frac{1}{(2\pi)^n} \left(\ln \|x\|_K\right)^\wedge(\xi) d\xi.$$
Since $\mu$ is a probability measure, one can see that
\begin{equation}\label{probab  measure}
\int_{S^{n-1}}(\ln\|x\|_K)^\wedge (\theta) d\theta=-(2\pi)^n
\end{equation} for any infinitely smooth origin-symmetric star body $K$ (see \cite[Remark 3.2]{KKYY}).

In our next Lemma we prove that a representation similar to (\ref{logrepr}) holds for all
infinitely smooth bodies, with  $\mu$ being a signed measure.
\begin{Lem}\label{Lem:log in R^n}
Let $K$ be an infinitely smooth origin-symmetric star body in $\mathbb{R}^n$, then
\begin{equation} \label{logrepr in R^n}
\ln \|x\|_K =-\frac{1 }{(2\pi)^n}\int_{S^{n-1}} \ln |(x, \xi ) |  \left(\ln
\|x\|_K\right)^\wedge(\xi) d\xi + C_K,
\end{equation}
where $C_K$ is the constant from (\ref{const}). \proof

Since the body $K$ is infinitely smooth, by Theorem \ref{Thm:n-1-deriv},  $\left(\ln
\|x\|_K\right)^\wedge(\xi)$ is a continuous homogeneous function of degree $-n$ on
$\mathbb{R}^n\setminus\{0\}$.

Let $\phi$ be an even test function supported outside of the origin, then

$$\left\langle \left(\int_{S^{n-1}} \ln |(x, \xi ) |  \left(\ln \|x\|_K\right)^\wedge(\xi)
d\xi\right)^\wedge , \phi \right\rangle$$ $$= \left\langle \int_{S^{n-1}} \ln |(x, \xi ) |
\left(\ln \|x\|_K\right)^\wedge(\xi) d\xi , \hat{\phi}(x) \right\rangle$$
$$=\int_{\mathbb{R}^n}\left[ \int_{S^{n-1}} \ln |(x, \xi ) |  \left(\ln
\|x\|_K\right)^\wedge(\xi) d\xi \right]\hat{\phi}(x) dx$$ $$= \int_{S^{n-1}}
\left[\int_{\mathbb{R}^n} \ln |(x, \xi ) | \hat{\phi}(x) dx \right] \left(\ln
\|x\|_K\right)^\wedge(\xi) d\xi$$
%&&\qquad\qquad=\int_{S^{n-1}} \left[\int_{\mathbb{R}}  \ln |t | \int_{(x, \xi )=t} \hat{\phi}(x) dx
%dt \right] \left(\ln \|x\|\right)^\wedge(\xi) d\xi
%\end{eqnarray*}

Now compute the inner integral using Fubini's theorem and the connection between the Radon and
Fourier transforms:
%\begin{eqnarray*}
%\int_{\mathbb{R}^n}  \ln |(x, \xi ) | \hat{\phi}(x) dx&=&\int_{\mathbb{R}}  \ln |t | \int_{(x, \xi
%)=t} \hat{\phi}(x) dx dt\\
%&=&\frac{1}{2\pi}\int_{\mathbb{R}} ( \ln |t | )^\wedge (z) \left(\int_{(x, \xi )=t} \hat{\phi}(x)
%dx\right)^\wedge(z)dz \\
%&=&-\frac{1}{2}\int_{\mathbb{R}} |z|^{-1}  \hat{\hat{\phi}}(z\xi) dz \\
%&=&-{2}^{n-1}\pi^n \int_{\mathbb{R}} |z|^{-1}  {{\phi}}(z\xi) dz\\
%&=&-(2\pi)^n \int_{0}^\infty z^{-1}  {{\phi}}(z\xi) dz\\
%\end{eqnarray*}

$$\int_{\mathbb{R}^n}  \ln |(x, \xi ) | \hat{\phi}(x) dx=\int_{\mathbb{R}}  \ln |t | \int_{(x, \xi
)=t} \hat{\phi}(x) dx dt$$ $$=\frac{1}{2\pi}\int_{\mathbb{R}} ( \ln |t | )^\wedge (z)
\left(\int_{(x, \xi )=t} \hat{\phi}(x) dx\right)^\wedge(z)dz =-\frac{1}{2}\int_{\mathbb{R}}
|z|^{-1}  \hat{\hat{\phi}}(z\xi) dz $$
$$=-{2}^{n-1}\pi^n \int_{\mathbb{R}} |z|^{-1}  {{\phi}}(z\xi) dz=-(2\pi)^n \int_{0}^\infty z^{-1}  {{\phi}}(z\xi)
dz$$
%\end{eqnarray*}

Here we used the formula for the Fourier transform of $\ln|t|$ (see \cite[p.362]{GS})
\begin{eqnarray}\label{FT:ln|t|}
\left(\ln|z|\right)^\wedge(t)= -\pi |t|^{-1}
\end{eqnarray}
outside of the origin. Therefore, passing from polar  to Euclidean coordinates and recalling from
Theorem \ref{Thm:n-1-deriv}, that $\left(\ln \|x\|_K\right)^\wedge$ is a homogeneous function of
degree $-n$ on $\mathbb{R}^n\setminus \{0\}$, we get
%\begin{eqnarray*}
%&&\langle \left(\int_{S^{n-1}} \ln |(x, \xi ) |  \left(\ln \|x\|\right)^\wedge(\xi)
%d\xi\right)^\wedge , \phi \rangle\\
%&&\qquad\qquad= -(2\pi)^n  \int_{S^{n-1}} \left[\int_{0}^\infty z^{-1}  {{\phi}}(z\xi) dz\right] \left(\ln \|x\|\right)^\wedge(\xi) d\xi\\
%&&\qquad\qquad= -(2\pi)^n  \int_{\mathbb{R}^n}   {{\phi}}(y)  \left(\ln \|x\|\right)^\wedge(y) dy\\
%&&\qquad\qquad= -(2\pi)^n \langle \left(\ln \|x\|\right)^\wedge, \phi \rangle.
%\end{eqnarray*}

$$\langle \left(\int_{S^{n-1}} \ln |(x, \xi ) |  \left(\ln \|x\|_K\right)^\wedge(\xi)
d\xi\right)^\wedge , \phi \rangle $$ $$= -(2\pi)^n  \int_{S^{n-1}} \left[\int_{0}^\infty z^{-1}
{{\phi}}(z\xi) dz\right] \left(\ln \|x\|_K\right)^\wedge(\xi) d\xi$$
$$= -(2\pi)^n  \int_{\mathbb{R}^n}   {{\phi}}(y)  \left(\ln \|x\|_K\right)^\wedge(y) dy= -(2\pi)^n \langle \left(\ln \|x\|_K\right)^\wedge, \phi \rangle.
$$

It follows that $$ \left(\int_{S^{n-1}} \ln |(x, \xi ) |  \left(\ln \|x\|_K\right)^\wedge(\xi)
d\xi\right)^\wedge =-(2\pi)^n  \left(\ln \|x\|_K\right)^\wedge$$
 as distributions outside of the origin. Hence, the functions $-(2\pi)^n \ln \|x\|_K$  and
  $\int_{S^{n-1}} \ln |(x, \xi ) |  \left(\ln \|x\|_K\right)^\wedge(\xi) d\xi$
 may differ only by a polynomial. But $$\frac{1}{(2\pi)^n}\int_{S^{n-1}} \ln |(x, \xi ) |
\left(\ln \|x\|_K\right)^\wedge(\xi) d\xi+ \ln \|x\|_K$$ is a homogeneous function of degree zero,
therefore  this polynomial is some constant $C$, which is exactly the constant from Definition
\ref{Def:Lzero}, as computed in \cite{KKYY}.

\qed
\end{Lem}

%Recall Parseval's formula on the sphere (A.Koldobsky, \cite{K2}). If $K$ and $L$ are infinitely
%smooth bodies in $\mathbb{R}^n$, then
%$$(2\pi)^n \int_{S^{n-1}}\|\theta\|_K^{-p}\|\theta\|_L^{-n+p}d\theta=
%\int_{S^{n-1}}(\|\theta\|_K^{-p})^\wedge(\xi)(\|\theta\|_L^{-n+p})^\wedge(\xi) d\xi,$$ for
%$p\in(0,n)$. Note that because the bodies are infinitely smooth, the Fourier transforms in the
%right-hand side are continuous functions on the sphere.

Now we need a version of Parseval's formula for $L_0$. How does the formula of Lemma
\ref{Lem:Parseval} look like if we pass to the limit as $p\to 0$? The answer to this question is
given in our next Lemma. Even though in the proof we are using an argument based on Lemma
\ref{Lem:log in R^n}, one can obtain the following Lemma by taking the limit in Parseval's formula.

\begin{Lem}\label{0-Parseval}
Let $K$ and $L$ be infinitely smooth origin-symmetric star bodies in $\mathbb{R}^n$, then
\begin{eqnarray*}
- \frac{1}{(2\pi)^n} \int_{S^{n-1}}\left[ \int_L \ln|( x,\xi) | dx \right]\,
(\ln\|x\|_K)^\wedge(\xi)d\xi=\int_L(\ln\|x\|_K-C_K) dx
\end{eqnarray*}
\end{Lem}
\proof By Lemma \ref{Lem:log in R^n} we have
\begin{equation*}% \label{logrepr in R^n}
-\frac{1 }{(2\pi)^n}\int_{S^{n-1}} \ln |(x, \xi ) |  \left(\ln \|x\|_K\right)^\wedge(\xi) d\xi =\ln
\|x\|_K - C_K.
\end{equation*}
Integrating this equality over the body $L$ we get the statement of the Lemma.

 \qed

Now we prove the main result of this section.

\begin{Thm}
Let $K$ and $L$ be two origin-symmetric star bodies in $\mathbb{R}^n$ such that $(\mathbb{R}^n,
\|\cdot\|_K)$ embeds in $L_0$ and
\begin{eqnarray}\label{ln-inequality}
%\frac{\int_L\ln|( x,\xi) | dx }{\mathrm{vol}(L)}\le \frac{\int_K\ln|( x,\xi)| dx }{\mathrm{vol}(K)}
\Gamma^*_0 K\subset \Gamma^*_0 L
\end{eqnarray}
for every $\xi\in S^{n-1}$, then $$\mathrm{vol}(L)\le \mathrm{vol}(K).$$
\end{Thm}
\proof

%Using the assumptions of the theorem and the fact that all 3-dimensional normed spaces embed in
%$L_0$ we get:

Since $(\mathbb{R}^n, \|\cdot\|_K)$ embeds in $L_0$, there exist a probability measure $\mu_K$ on
$S^{n-1}$ (which is the restriction of the Fourier transform of $\ln\|x\|_K$ to the unit sphere)
and a constant $C_K$ from Definition \ref{Def:Lzero}.

Rewrite   inequality (\ref{ln-inequality}) as follows:
\begin{eqnarray*}
\frac{\int_L\ln|( x,\xi) | dx }{\mathrm{vol}(L)}\le \frac{\int_K\ln|( x,\xi)| dx
}{\mathrm{vol}(K)},
%\Gamma^*_0 K\subset \Gamma^*_0 L
\end{eqnarray*}
and integrate it over $S^{n-1}$ with respect to $\mu_K$ to get
\begin{eqnarray*}
\int_{S^{n-1}} \frac{\int_L\ln|( x,\xi) | dx }{\mathrm{vol}(L)}d\mu_K(\xi)\le \int_{S^{n-1}}
\frac{\int_K\ln|( x,\xi) | dx }{\mathrm{vol}(K)}d\mu_K(\xi)
\end{eqnarray*}
Using  the Fubini theorem and the definition of embedding in $L_0$, %Lemma \ref{0-Parseval},
we get
$$\frac{1}{\mathrm{vol}(L)}\int_L(\ln\|x\|_K-C_K) dx \le \frac{1}{\mathrm{vol}(K)}\int_K(\ln\|x\|_K-C_K)
dx$$

Therefore
$$\frac{1}{\mathrm{vol}(L)}\int_L\ln\|x\|_K dx \le \frac{1}{\mathrm{vol}(K)}\int_K \ln\|x\|_K
dx=-\frac{1}{n},$$ where the latter equality follows from the formula $$\frac{1}{\mathrm{vol}(K)}
\int_K \|x\|_K^p d x = \frac{n}{n+p},$$ that we had earlier, after differentiating and letting
$p=0$.

Now use the following inequality from Milman and Pajor \cite[Section 2.2]{MP}:
\begin{equation}\label{MP-log}
 \frac{1}{\mathrm{vol}(L)}\int_L\ln\|x\|_K dx\ge -\frac{1}{n} + \frac{1}{n} [\ln(\mathrm{vol}(L)) -
 \ln(\mathrm{vol}(K))]
 \end{equation}
Therefore
$$\mathrm{vol}(L) \le \mathrm{vol}(K).$$

\qed

\begin{Cor}
Since every three dimensional normed space embeds in $L_0$ (see \cite[Corollary 4.3]{KKYY}), the
previous theorem holds for all convex bodies in $\mathbb{R}^3$.
\end{Cor}

To prove our next Theorem we need the following Lemma.

\begin{Lem}\label{FT: norm to -n} Let $K$ be an origin-symmetric star body in $\mathbb{R}^n$, then
the Fourier transform of $\|x\|_K^{-n}$ is a continuous function on  $\mathbb{R}^n\setminus \{0\}$
and equals
\begin{eqnarray*}
(\|x\|_K^{-n})^\wedge(\xi) = &-& n \int_K\ln|( x,\xi) | dx+\\
&+& (n \Gamma'(1)-1)\mathrm{vol}(K) -\int_{S^{n-1}}\|\theta\|_K^{-n} \ln \|\theta\|_K  d\theta.
\end{eqnarray*}
\end{Lem}
\proof Let $\phi$ be an even test function. Using the definition of the action of a homogeneous
function of degree $-n$ (see \cite[p.303]{GS}) we get
\begin{eqnarray*}
&&\hspace{-1cm}\langle (\|x\|_K^{-n})^\wedge,\phi\rangle   = \langle \|x\|_K^{-n},\hat{\phi}(x)\rangle\\
                                            &=& \int_{B_1(0)} \|x\|_K^{-n}(\hat{\phi}(x)-\hat{\phi}(0)) dx+\int_{\mathbb{R}^n\setminus B_1(0)}
                                            \|x\|_K^{-n}\hat{\phi}(x)dx\\
                                            &=& \int_{S^{n-1}} \int_0^1  r^{-1}
                                            \|\theta\|_K^{-n}(\hat{\phi}(r\theta)-\hat{\phi}(0))dr
                                            d\theta + \int_{S^{n-1}} \int_1^\infty  r^{-1}
                                            \|\theta\|_K^{-n}\hat{\phi}(r\theta)dr
                                            d\theta\\
                                            &=& \int_{S^{n-1}} \|\theta\|_K^{-n}\left( \int_0^1  r^{-1}
                                            (\hat{\phi}(r\theta)-\hat{\phi}(0))dr +  \int_1^\infty  r^{-1}
                                            \hat{\phi}(r\theta)dr\right)
                                            d\theta\\
                                            &=&\frac12 \int_{S^{n-1}} \|\theta\|_K^{-n} \langle |r|^{-1},
                                            \hat{\phi}(r\theta)\rangle d\theta\\
                                            &=&\frac12 \int_{S^{n-1}} \|\theta\|_K^{-n} \langle
                                            2\Gamma'(1) - 2\ln |t|, \int_{(\theta,\xi)=t}
                                            \phi(\xi)d\xi\rangle d\theta\\
                                            &=&\langle \int_{S^{n-1}} \|\theta\|_K^{-n}
                                            \left(\Gamma'(1) - \ln |(\theta,\xi)| \right)d\theta,
                                            \phi(\xi)\rangle ,
\end{eqnarray*}
here we used the formula for the Fourier transform of $|r|^{-1}$ from \cite[p.361]{GS}:
$$(|r|^{-1})^\wedge(t)= 2\Gamma'(1) - 2\ln |t|.$$
Thus we have proved that
\begin{eqnarray}\label{FT -n} (\|x\|_K^{-n})^\wedge(\xi) &=& \int_{S^{n-1}}
\|\theta\|_K^{-n}\Big(\Gamma'(1)-\ln |(\theta,\xi)|\Big)d\theta.
\end{eqnarray}

Next,  let us compute the following:
\begin{eqnarray*}
&&\hspace{-1cm}\int_K\ln|( x,\xi) | dx=  \int_{S^{n-1}} \int_0^{\|\theta\|_K^{-1}} r^{n-1}
\ln|(r\theta,\xi)| dr d\theta\\
&=& \int_{S^{n-1}} \int_0^{\|\theta\|_K^{-1}}  r^{n-1}\ln r dr d\theta + \int_{S^{n-1}}
\ln|(\theta,\xi)| \int_0^{\|\theta\|_K^{-1}} r^{n-1}  dr d\theta\\
 &=&- \frac1n\int_{S^{n-1}}\Big( \|\theta\|_K^{-n} \ln \|\theta\|_K + \frac1n  \|\theta\|_K^{-n} \Big) d\theta +
 \frac1n \int_{S^{n-1}} \|\theta\|_K^{-n}\ln|(\theta,\xi) | d\theta
\end{eqnarray*}
Therefore
\begin{eqnarray*}
&&\hspace{-1.5cm}\int_{S^{n-1}} \|\theta\|_K^{-n}\ln|(\theta,\xi) | d\theta=\\
&&=n \int_K\ln|( x,\xi) | dx+\int_{S^{n-1}} \Big( \|\theta\|_K^{-n} \ln \|\theta\|_K + \frac1n
\|\theta\|_K^{-n} \Big) d\theta.
\end{eqnarray*}
Combining this formula with the formula (\ref{FT -n}), we get
\begin{eqnarray*}
(\|x\|_K^{-n})^\wedge(\xi) = &-& n \int_K\ln|( x,\xi) | dx+\\
&+& (n \Gamma'(1)-1)\mathrm{vol}(K) -\int_{S^{n-1}}\|\theta\|_K^{-n} \ln \|\theta\|_K  d\theta.
\end{eqnarray*}

\qed

\begin{Thm}
There are convex bodies $K$ and $L$  in $\mathbb{R}^n$, $n\ge 4$ such that
\begin{eqnarray*}
%\frac{\int_L\ln|( x,\xi) | dx }{\mathrm{vol}(L)}\le \frac{\int_K\ln|( x,\xi)| dx }{\mathrm{vol}(K)}
\Gamma^*_0 K\subset \Gamma^*_0 L
\end{eqnarray*}
%\begin{eqnarray*}
%\frac{\int_K\ln|( x,\xi) | dx }{\mathrm{vol}(K)}\le \frac{\int_L\ln|( x,\xi) | dx
%}{\mathrm{vol}(L)}
%\end{eqnarray*}
for every $\xi\in S^{n-1}$, but  $$\mathrm{vol}(K)< \mathrm{vol}(L).$$
\end{Thm}
\noindent{\bf Proof.}
 Let $L$ be   a strictly convex infinitely smooth body in $\mathbb{R}^n$, $n\ge 4$, for which $-(\ln
\|x\|_L)^\wedge$ is not positive everywhere. (See  \cite[Theorem 4.4]{KKYY} for an explicit
construction of such a body).

Let $\xi\in S^{n-1}$ be such that $-(\ln \|x\|_L)^\wedge(\xi)<0$. By continuity of the function
$(\ln\|x\|_L)^\wedge(\theta)$ on the sphere there is a neighborhood  of $\xi$ where this function
is negative. Let
$$\Omega=\{\theta \in S^{n-1}: -(\ln \|x\|_L^{})^\wedge(\theta)<0\}.$$
Choose an infinitely smooth body $D$ whose Minkowski norm $\|x\|_D$ is equal to 1 outside of
$\Omega$ and $\|x\|_D< 1 $ for $x\in \Omega$. Let $v$ be a homogeneous function of degree $0$ on
$\mathbb{R}^n\setminus\{0\}$, defined as follows:
$$v(x)=\ln\|x\|_D-\ln|x|_2.$$
Clearly $v(x)< 0$ if $x\in \Omega$ and $v(x)= 0$ if $x\in S^{n-1}\setminus \Omega$.

In view of   Theorem \ref{Thm:n-1-deriv}, the Fourier transforms of $\ln\|x\|_D$ and $\ln|x|_2$
outside of the origin are some homogeneous functions of degree $-n$, therefore  the Fourier
transform of $ v(x)$ outside of the origin is equal to $|x|_2^{-n}g(x/|x|_2)$ for some infinitely
smooth function $g$ on $S^{n-1}$. Since by (\ref{probab  measure})
$$\int_{S^{n-1}} (\ln\|x\|_D)^\wedge (\theta) d\theta=\int_{S^{n-1}} (\ln|x|_2)^\wedge (\theta) d\theta=-(2\pi)^n,$$
we have
\begin{equation}\label{int of g}
\int_{S^{n-1}} g(\theta) d\theta=0.
\end{equation}

Define a body $K$ by the formula:
\begin{eqnarray}\label{def2 K}
\frac{\|x\|_K^{-n}}{\mathrm{vol}(K)}=\frac{\|x\|_L^{-n}}{\mathrm{vol}(L)}  + n (2\pi)^{-n} \epsilon
|x|_2^{-n}g({x}/{|x|_2}).
\end{eqnarray}
%Note that, because of the formula (\ref{int of g}) $\mathrm{vol}(K)=c$.
Note that formula (\ref{int of g}) validates this definition, since integrating the last equality
over the unit sphere we get the same quantity in both sides. Also, since $L$ is strictly convex,
there is an $\epsilon$ small enough, so that $K$ is also convex (see e.g. the perturbation argument
from \cite[Section 5.1]{K3}). From now on we fix such an $\epsilon$.

Now we will show that $K$ together with $L$ constructed above satisfy the assumptions of the
theorem. Apply the Fourier transform to both sides of (\ref{def2 K}). Note, that  the Fourier
transform of $|x|_2^{-n}g({x}/{|x|_2})$ is  equal to $(2\pi)^n v$ on test functions, whose Fourier
transform  is supported outside of the origin. Such distributions can differ only by a polynomial,
which must be a constant in this case, since both functions cannot grow faster than a logarithm
(see Lemma \ref{FT: norm to -n}). So
$$\left(|x|_2^{-n}g({x}/{|x|_2})\right)^\wedge=(2\pi)^n (v+ \alpha),$$
for some constant $\alpha$ whose value has no significance for us. Hence, by Lemma \ref{FT: norm to
-n}, the Fourier transform of (\ref{def2 K}) looks as follows:

\begin{eqnarray}\label{with C}
-\frac{n\int_K\ln|( x,\xi) | dx }{\mathrm{vol}(K)}=- \frac{n\int_L\ln|( x,\xi) | dx
}{\mathrm{vol}(L)} +n  \epsilon\cdot v(\xi) + C,
\end{eqnarray}
where the constant $C$ equals
\begin{eqnarray*}
C=\frac{\int_K \|\theta\|_K^{-n} \ln \|\theta\|_K d\theta dx }{\mathrm{vol}(K)}- \frac{\int_L
\|\theta\|_L^{-n} \ln \|\theta\|_L d\theta dx }{\mathrm{vol}(L)} + n  \epsilon\cdot \alpha.
\end{eqnarray*}

Since the bodies $L$ and $D$ are fixed, dilating the body $K$ we can make this constant equal to
zero. Indeed, multiply the Minkowski functional of $K$ by a positive constant $\lambda$, then
\begin{eqnarray*}
C&=&\frac{\int_K (\lambda\|\theta\|_K) ^{-n} \ln \lambda \|\theta\|_K d\theta dx }{\lambda^{-n}
\mathrm{vol}(K)}- \frac{\int_L \|\theta\|_L^{-n} \ln \|\theta\|_L d\theta dx }{\mathrm{vol}(L)} +
n  \epsilon\cdot\alpha\\
&=&\frac{\int_K \|\theta\|_K ^{-n} \left[\ln \lambda+\ln  \|\theta\|_K \right]d\theta dx }{
\mathrm{vol}(K)}- \frac{\int_L \|\theta\|_L^{-n} \ln \|\theta\|_L d\theta dx }{\mathrm{vol}(L)} +
n  \epsilon\cdot\alpha\\
&=&n \ln \lambda +\frac{\int_K \|\theta\|_K ^{-n} \ln  \|\theta\|_K d\theta dx }{ \mathrm{vol}(K)}
- \frac{\int_L \|\theta\|_L^{-n} \ln \|\theta\|_L d\theta dx }{\mathrm{vol}(L)} + n
\epsilon\cdot\alpha
\end{eqnarray*}
One can choose a $\lambda>0$ so that $C=0$. Therefore from (\ref{with C}) we get
\begin{eqnarray}\label{eqn:inequality for sections}
\frac{\int_K\ln|\langle x,\xi\rangle | dx }{\mathrm{vol}(K)}=\frac{\int_L\ln|\langle x,\xi\rangle |
dx }{\mathrm{vol}(L)}  - \epsilon\ v(\xi)\ge \frac{\int_L\ln|\langle x,\xi\rangle | dx
}{\mathrm{vol}(L)},
\end{eqnarray}
since $v$ is non-positive. Therefore
\begin{eqnarray*}
%\frac{\int_L\ln|( x,\xi) | dx }{\mathrm{vol}(L)}\le \frac{\int_K\ln|( x,\xi)| dx }{\mathrm{vol}(K)}
\Gamma^*_0 K\subset \Gamma^*_0 L.
\end{eqnarray*}

Now using  Parseval's formula  and  inequality (\ref{eqn:inequality for sections}) we get
\begin{eqnarray*}
&&\hspace{-1.5cm}\frac{1}{\mathrm{vol}(K)}\int_K(\ln\|x\|_L-C_L) dx=\\
&=&-\frac{1}{(2\pi)^n}\frac{1}{\mathrm{vol}(K)} \int_{S^{n-1}} \left[
\int_K \ln|\langle x,\xi\rangle | dx\right]  (\ln\|x\|_L)^\wedge(\xi)d\xi\\
&=& -\frac{1}{(2\pi)^n}\int_{S^{n-1}} \left[\frac{1}{\mathrm{vol}(L)}\int_L \ln|\langle
x,\xi\rangle | dx - \epsilon v(\xi) \right]\, (\ln\|x\|_L)^\wedge(\xi)d\xi\\
&=&-\frac{1}{(2\pi)^n}\frac{1}{\mathrm{vol}(L)}\int_{S^{n-1}} \left[\int_L \ln|\langle x,\xi\rangle
| dx \right]\, (\ln\|x\|_L)^\wedge(\xi)d\xi\\
&&+\frac{1}{(2\pi)^n}\frac{1}{\mathrm{vol}(L)}\int_{S^{n-1}}\epsilon v(\xi) (\ln\|x\|_L)^\wedge(\xi)d\xi\\
&<&-\frac{1}{(2\pi)^n}\frac{1}{\mathrm{vol}(L)}\int_{S^{n-1}} \left[\int_L \ln|\langle x,\xi\rangle
| dx \right]\,
(\ln\|x\|_L)^\wedge(\xi)d\xi\\
&=&\frac{1}{\mathrm{vol}(L)}\int_L(\ln\|x\|_L-C_L) dx,
\end{eqnarray*}
where the inequality follows from the fact that $v$ is non-positive and it is supported on the set
where $-(\ln\|x\|_L)^\wedge(\xi)<0$.

Recalling the inequality (\ref{MP-log})
 $$-\frac1n\ge\frac{1}{\mathrm{vol}(K)}\int_K\ln\|x\|_L dx\ge -\frac{1}{n} + \frac{1}{n} [\ln(\mathrm{vol}(K)) -
 \ln(\mathrm{vol}(L))],$$
we get
$$\mathrm{vol}(K) < \mathrm{vol}(L).$$

\qed

{\bf Acknowledgments.} The authors wish to thank A.~Koldobsky for reading this manuscript and
making many valuable suggestions.

\end{document}